\documentclass[a4paper,12pt,centertags,psamsfonts]{amsart}
\usepackage{amsmath,latexsym}
\usepackage{amssymb}
\usepackage{times}
\DeclareSymbolFont{SY}{U}{psy}{m}{n}
\DeclareMathSymbol{\emptyset}{\mathord}{SY}{'306}
\chardef\bslash=`\\ 

\newtheorem*{thm}{Theorem} 

\newtheorem{prop}{Proposition}

\theoremstyle{definition}
\newtheorem{defn}{Definition}[]

\theoremstyle{remark}
\newtheorem*{rem}{Remark} 
\newtheorem*{notation}{Notation}

\newcommand{\propref}[1]{Proposition~\ref{#1}}

\newcommand{\Na}{\mathbb{N}}

\newcommand{\Zi}{\mathbb{Z}}

\newcommand{\R}{\mathbb{R}}
\newcommand{\fR}{\mathfrak{R}}

\newcommand{\C}{\mathbb{C}}

\let\bsy\boldsymbol
\def\P{\bsy P}
\def\T{\bsy T}
\def\t{\bsy t}

\def\p{\bsy p}
\def\f{\bsy f}
\def\K{\bsy K}
\def\k{\bsy k}

\def\F{\bsy F}

\def\H{\mathcal{H}}

\newcommand{\Ss}{\mathcal{S}}

\DeclareMathOperator{\sgn}{sign}

\newcommand{\eval}[2][\right]{\relax
  \ifx#1\right\relax \left.\fi#2#1\rvert}

\let\abs=\envert

\let\norm=\enVert

\let\hnorm=\henVert
\newcommand{\kskob}[1]{\left[#1\right]}
\let\ksk=\kskob

\let\sk=\skob
\newcommand{\skoba}[1]{\left\langle#1\right\rangle}
\let\ska=\skoba

\let\hska=\hskoba

\let\sets=\fskoba
\begin{document}
\renewcommand{\sectionmark}[1]{}
\title[Simultaneous unitary equivalence to Carleman
operators]{Simultaneous unitary equivalence to Carleman
operators with arbitrarily smooth kernels}
\author[I. M. Novitski\u i]{Igor M. Novitski\u i}
\address{Institute for Applied Mathematics, Russian Academy of Sciences,
92, Zaparina Street, Khabarovsk 680 000, Russia}
\email{novim@iam.khv.ru}
\thanks{Research supported in part by grant N 03-1-0-01-009 from
the Far-Eastern Branch of the Russian Academy of Sciences. This paper
was written in November 2003, when the author enjoyed the hospitality of the
Mathematical Institute of Friedrich-Schiller-University, Jena, Germany}
\begin{abstract}
In this paper, we describe families of those bounded linear operators
on a separable Hilbert space that are simultaneously unitarily equivalent to
integral operators on $L_2(\R)$ with bounded and \textit{arbitrarily} smooth
Carleman kernels. The main result is a qualitative sharpening of an earlier
result of \cite{nov:Isra}.
\end{abstract}
\maketitle

\section{Introduction. Main result}
Throughout,  $\H$ will denote a separable Hilbert space with the inner product
$\hska{\cdot,\cdot}$ and the norm $\hnorm{\cdot}$, $\fR(\H)$ the algebra
of all bounded linear operators on $\H$, and $\C$, and $\Na$, and $\Zi$,
the complex plane, the set of all positive integers, the set of all integers,
respectively. For an operator $A$ in $\fR(\H)$, $A^*$ will denote the
Hilbert space adjoint of $A$ in $\fR(\H)$.

Throughout,  $C(X,B)$, where $B$ is a Banach space (with norm
$\norm{\cdot}_B$), denote the Banach space (with the norm $\norm{f}_{C(X,B)}
=\sup\limits_{x\in X}\,\norm{f(x)}_B$) of
continuous $B$-valued functions defined on a
locally compact space $X$ and \textit{vanishing at infinity\/} (that is, given
any $f \in C(X,B)$ and $\varepsilon>0$, there exists a compact subset
$X(\varepsilon,f) \subset X$ such that $\norm{f(x)}_{B} < \varepsilon$
whenever $x \not\in  X(\varepsilon,f)$).

Let $\R$ be the real line $(-\infty,+\infty)$ with the Lebesgue measure,
and let $L_2=L_2(\R)$ be the Hilbert
space of (equivalence classes of) measurable complex-valued functions on
$\R$ equipped with the inner product
$$
\ska{f,g}=\int_\R f(s)\overline{g(s)}\,ds
$$
and the norm
$\norm{f}=\ska{f,f}^{\frac{1}2}$.

A linear operator
                 $T : L_2 \to  L_2$
is said to be \textit{integral\/}
if there exists a measurable function
$\T$ on the Cartesian product $\R^2=\R\times\R$, a \textit{kernel\/},
such that, for every $f\in L_2$,
$$
               (Tf)(s)=\int_{\R} \T(s,t)f(t)\,dt
$$
for almost every $s$ in $\R$. A kernel $\T$ on
$\R^2$ is said to be \textit{Carleman\/} if $\T(s,\cdot) \in L_2$
for almost every fixed $s$ in $\R$.
An integral operator with a kernel $\T$ is
called \textit{Carleman\/} if $\T$ is a Carleman kernel.
Every Carleman kernel, $\T$, induces a \textit{Carleman
function\/} $\t$ from $\R$ to $L_2$ by
$\t(s)=\overline{\T(s,\cdot)}$
for all $s$ in ${\R}$ for which $\T(s,\cdot)\in L_2$.

The integral representability problem for linear operators
stems from the work \cite{Neu} of von~Neumann, and
is now well enough understood. The problem involves the question: which operators
are unitarily equivalent to an integral operator?
Now we recall a characterization of Carleman representable operators
to within unitary equivalence \cite[p.~99]{Kor:book1}, \cite[Section 15]{Halmos:Sun}:
\begin{prop}\label{Kor}
A necessary and sufficient condition that
an operator $S\in \fR(\H)$ be unitarily equivalent to an integral operator
with Carleman kernel
is that there exist an
orthonormal sequence $\left\{e_n\right\}$ such that
$$
\hnorm{S^*e_n}\rightarrow 0\quad\text{as $n\rightarrow\infty$}
$$
(or, equivalently, that $0$ belong to the right essential spectrum of $S$).
\end{prop}

Given any non-negative integer $m$,
we impose on a Carleman kernel $\K$ the following smoothness conditions:
\begin{enumerate}
\renewcommand{\labelenumi}{(\roman{enumi})}
\item the function $\K$ and all its partial derivatives
on $\R^2$ up to order $m$ are in $C(\R^2,\C)$,
\item the Carleman function $\k$,
$\k(s)=\overline{\K(s,\cdot)}$,
and all its (strong) derivatives on ${\R}$ up to order $m$
are in $C(\R,L_2)$.
\end{enumerate}

\begin{defn} A function $\K$ that satisfies Conditions (i), (ii)
is called a \textit{$SK^m$-kernel\/} \cite{nov:Isra}.
\end{defn}

Now we are  in a position to formulate our result on
simultaneous integral representability of operator families
by $SK^m$-kernels.

\begin{prop}[\cite{nov:Isra}]\label{msmooth}
If for a countable family
$\sets{B_r\mid r\in\Na}\subset \fR(\H)$
there exists an orthonormal sequence
$\sets{e_n}$
such that
$$
\sup_{r\in\Na}\hnorm{B^*_r e_{n}}\rightarrow0\quad
 \text{as\ } n \rightarrow \infty,
$$
then for each fixed non-negative integer $m$ there exists
a unitary operator $U_m:\H\to L_2$ such that all the operators
$U_mB_r U_m^{-1}$ $(r\in\Na)$ are bounded Carleman operators having
$SK^m$-kernels.
\end{prop}

In \cite{nov:Isra}, there is a counterexample which shows that
\propref{msmooth} may fail to be true if the family
$\sets{B_r}$ is not countable.

The purpose of this paper is to restrict the conclusion of \propref{msmooth}
to \textit{arbitrarily\/} smooth Carleman kernels. Now define these kernels.

\begin{defn}
We say  that a function $\K$ is a \textit{$SK^\infty$-kernel}
(\cite{Je}, \cite{SP})
if it is a $SK^m$-kernel for each non-negative integer $m$.
\end{defn}

\begin{thm}\label{arsmooth}
If for a countable family
$\sets{B_r\mid r\in\Na}\subset \fR(\H)$
there exists an orthonormal sequence
$\sets{v_n}$
such that
\begin{equation}\label{main}
\sup_{r\in\Na}\hnorm{B^*_r v_n}\to0\quad
 \text{as $n\to\infty$},
\end{equation}
then there exists
a unitary operator $U_\infty:\H\to L_2$ such that all the operators
$U_\infty B_r U_\infty^{-1}$ $(r\in\Na)$ are Carleman operators having
$SK^\infty$-kernels.
\end{thm}

This theorem, which is our main result, will be proved in the next section
of the present paper. The proof yields an explicit construction of the unitary operator
$U_\infty:\H\to L_2$.
The construction of $U_\infty$ is
independent of those spectral points of $B_r$ ($r\in\Na$) that are
different from $0$, and is defined by $U_\infty f_n=u_n$ ($n\in\Na$),
where $\{f_n\}$, $\{u_n\}$ are orthonormal bases in $\H$ and $L_2$, respectively,
whose elements can be explicitly described in terms of the operator family.

\section{Proof of Theorem}
The proof has two steps.
\subsection*{Step 1.}
Assume that
\begin{equation*}\label{norm}
\sup\limits_{r\in\Na}\norm{B_r}\le1.
\end{equation*}
This is a harmless
assumption, involving no loss of generality; just replace
$B_r$ with $\norm{B_r}>1$ by $\dfrac{B_r}{\norm{B_r}}$.
Find a subsequence
$\sets{e_k}_{k=1}^\infty$
of the sequence
$\sets{v_n}$ in \eqref{main}
so that
\begin{equation}\label{null-seq}
\begin{gathered}
\sum_k\sup_{r\in\Na}\hnorm{S_r^*e_k}^{\frac14}\le
\sum_k\sup_{r\in\Na}\hnorm{rS_r^*e_k}^{\frac14}\\=
\sum_k\sup_{r\in\Na}\hnorm{B_r^*e_k}^{\frac14}=M<\infty,
\end{gathered}
\end{equation}
where $S_r=\dfrac1rB_r$ ($r\in\Na$)
(the sum notation $\sum\limits_k$ will always be used instead of
the more detailed symbol $\sum\limits_{k=1}^\infty$).
For each $r$, let
\begin{equation}\label{E}
Q_r=(1-E)S_r,\quad J_r=S_r^*E,
\end{equation}
where
$E$ is the orthogonal projection onto
the closed linear span $H$ of the $e_k$'s,
and observe that
\begin{equation}\label{splitting}
S_r=Q_r+J^*_r.
\end{equation}
Assume, with no loss of generality, that
$\dim (1-E)H=\infty$, and let
$\{e_k^\perp\}_{k=1}^\infty$ be any orthonormal basis for $(1-E)H$.
Let $\{f_n\}_{n=1}^\infty$ denote any basis in $\H$
consisting of the elements of the set
$\sets{e_k}\cup\sets{e_k^\perp}$.
It follows from \eqref{null-seq} that
$$
\sum_n\hnorm{J_rf_n}=\sum_k\hnorm{ J_re_k}
\le\sum_k\sup_{r\in\Na}\hnorm{S_r^* e_k}\leq M^4,
$$
and hence
that $J_r$ and
$J_r^*$ are Hilbert--Schmidt operators, for each $r$.

For each $h\in\H$, let
\begin{equation}\label{defdknums}
d(h)=\sup_{r\in\Na}\hnorm{J_rh}^\frac{1}4+\sup_{r\in\Na}\hnorm{J_r^*h}^\frac{1}4+
\sup_{r\in\Na}\hnorm{\varGamma_rh},
\end{equation}
where, for each $r$,
\begin{equation}\label{lamgam}
\varGamma_r=\Lambda S_r, \text{\ and\ }
\Lambda =\sum_k\frac1k\hska{\cdot,e_k^\perp} e_k^\perp.
\end{equation}
It is clear that $\Lambda$ and $\varGamma_r$ ($r\in\Na$) are
Hilbert-Schmidt operators on $\H$.
 Prove that
\begin{equation}\label{to0}
 d(e_k)\to 0\quad\text{as $k\to\infty$}.
\end{equation}
Using known facts about Hilbert--Schmidt operators
(see \cite[Chapter III]{Goh:Kr}), write the following
relations
\begin{equation}\label{2.3}
\begin{gathered}
\sum_r\sup_{r\in\Na}\hnorm{J^*_r e_k}^2 \leq
\sum_r\sum_k\hnorm{J^*_r e_k}^2 \leq
\sum_r\bsy{\vert}J^*_r\bsy{\vert}_2^2
=\sum_r\bsy{\vert}J_r\bsy{\vert}_2^2\\
=\sum_r\sum_k\hnorm{J_r e_k}^2
=\sum_r\sum_k\hnorm{S^*_r e_k}^2\\
\le\sum_r\dfrac1{r^2}\sum_k\sup_{r\in\Na}\hnorm{rS^*_r e_k}^2
\le\dfrac{M^8\pi^2}6,
\end{gathered}
\end{equation}
where
$\bsy{\vert}\cdot\bsy{\vert}_2$ is the Hilbert--Schmidt norm.
Observe also that
\begin{equation}\label{2.4}
\begin{gathered}
\sum_k\sup_{r\in\Na}\hnorm{\varGamma_r e_k}^2
\le
\sum_r\sum_k\hnorm{\varGamma_r e_k}^2
\\
\leq
\sum_r\bsy{\vert}\varGamma_r\bsy{\vert}^2_2
=\sum_r\bsy{\vert}\varGamma_r^*\bsy{\vert}^2_2
=\sum_r\sum_n\hnorm{S^*_r\Lambda f_n}^2\\
\leq\sum_r\dfrac1{r^2}\sum_k\hnorm{\Lambda e_k^\perp}^2
=\sum_r\dfrac1{r^2}\sum_k\frac1{k^2}=\dfrac{\pi^4}{36}.
\end{gathered}
\end{equation}
Then \eqref{to0} follows immediately from \eqref{2.3}, \eqref{2.4},
\eqref{null-seq}, and \eqref{E}.

\begin{notation}
If an equivalence class
$f\in L_2$ contains a function belonging to $C(\R,\C)$, then we shall use
$\ksk{f}$ to denote that function.
\end{notation}
Take any orthonormal basis $\{u_n\}$ for $L_2$ which satisfies conditions:
\begin{enumerate}
\renewcommand{\labelenumi}{(\alph{enumi})}
\item the terms of the derivative sequence
$\left\{\ksk{u_n} ^{(i)}\right\}$ are in $C(\R,\C)$, for each $i$
(here and throughout, the letter $i$ is reserved for all non-negative integers),

\item $\{u_n\}=\{g_k\}_{k=1}^\infty\cup\{h_k\}_{k=1}^\infty$, where
$\{g_k\}_{k=1}^\infty\cap\{h_k\}_{k=1}^\infty=\varnothing$,
and,  for each $i$,
\begin{equation}\label{hki}
\sum_k H_{k,i}<\infty\quad
\text{with $H_{k,i}=\norm{\ksk{h_k}^{(i)}}_{C(\R,\C)}$}\quad (k\in\Na),
\end{equation}
\item there exist a subsequence
$\left\{x_k\right\}_{k=1}^\infty\subset\{e_k\}$ and
a strictly increasing sequence $\left\{n(k)\right\}_{k=1}^\infty$
of positive integers
such that, for each $i$,
\begin{gather}\label{zndn}
\sum_k d(x_k)\left(G_{k,i}+1\right)<\infty\quad
\text{with $G_{k,i}=\norm{\ksk{g_k}^{(i)}}_{C(\R,\C)}$}\quad (k\in\Na),\\
\label{sumrk}
\sum_k kH_{n(k),i}<\infty.
\end{gather}
\end{enumerate}

\begin{rem}
Let $\{ u_n\}$ be an orthonormal basis for $L_2$ such that, for
each $i$,
\begin{gather}\label{1}
\kskob{u_n}^{(i)}\in C(\R,\C)\quad(n\in\Na),\\
\label{2}\norm{\kskob{u_n}^{(i)}}_{C(\R,\C)}\le D_nA_i\quad(n\in\Na),\\
\label{3}
\sum_kD_{n_k}<\infty,
\end{gather}
where $\{D_n\}_{n=1}^\infty$, $\{A_i\}_{i=0}^\infty$ are sequences
of positive numbers, and $\{n_k\}_{k=1}^\infty$ is a subsequence
of $\Na$ such that $\Na\setminus\sets{n_k}_{k=1}^\infty$ is a countable set.
By \eqref{to0}, the basis $\{u_n\}$
satisfies Conditions (a)-(c) with $h_k=u_{n_k}$
($k\in\Na$) and
$\sets{g_k}_{k=1}^\infty=\{u_n\}\setminus\{h_k\}_{k=1}^\infty$.

To show the existence of
a basis $\{u_n\}$ satisfying \eqref{1}-\eqref{3}, consider
a Lemari\'e-Meyer wavelet,
$$
u(s)=\dfrac1{2\pi}\int_{\R}e^{i\xi(\frac12+s)}
\sgn\xi b(|\xi|)\,d\xi\quad (s\in\R),
$$
with the bell function $b$ belonging to
$C^\infty(\R)$ (for construction of the Lemari\'e-Meyer wavelets
we refer to \cite{LeMe}, \cite[\S~4]{Ausch}, \cite[Example D, p.~62]{Her}).
In this case, $u$ belongs to the Schwartz class $\Ss(\R)$, and hence
all the derivatives $\kskob{u}^{(i)}$
are in $C(\R,\C)$).
The ``mother function'' $u$ generates an orthonormal basis for $L_2$ by
$$
u_{jk}(s)=2^{\frac j2}u(2^js-k)\quad  (j,\,k\in\Zi).
$$
 Rearrange, in a completely arbitrary manner, the orthonormal set
$\{u_{jk}\}_{j,\,k\in\Zi}$ into a simple sequence,
so that it becomes $\{ u_n\}_{n\in\Na}$. Since, in view of this rearrangement,
to each $n\in\Na$ there corresponds a unique pair of integers
$j_n$, $k_n$, and
conversely, we can write, for each $i$,
$$
\norm{\kskob{u_n}^{(i)}}_{C(\R,\C)}=\norm{\kskob{u_{j_nk_n}}^{(i)}}_{C(\R,\C)}
\le D_nA_i,
$$
where
$$
D_n=\begin{cases}
2^{j_n^2}&\text{if $j_n>0$,}\\
\left(\dfrac1{\sqrt{2}}\right)^{\abs{j_n}}&\text{if $j_n\le0$,}
\end{cases}
\qquad A_i=2^{\left(i+\frac12\right)^2}\norm{\kskob{u}^{(i)}}_{C(\R,\C)}.
$$
Whence it follows that if $\{n_k\}_{k=1}^\infty\subset\Na$ is a subsequence
such that $j_{n_k}\to -\infty$ as $k\to\infty$, then
$$\sum_kD_{n_k}<\infty.$$
Thus,  the basis $\{u_n\}$ satisfies Conditions \eqref{1}-\eqref{3}.
\end{rem}

Let us return to the proof.
Let $\sets{x_k^\perp}_{k=1}^\infty=\sets{e_k^\perp}_{k=1}^\infty
\cup(\sets{e_k}_{k=1}^\infty\setminus\sets{x_k}_{k=1}^\infty)$, and observe
that $\sets{f_n}_{n=1}^\infty=\sets{x_k}_{k=1}^\infty\cup\sets{x_k^\perp}_{k=1}^\infty$.

Now construct a candidate for the desired
unitary operator in the theorem.
Define a unitary operator $U_\infty:\H\to L_2$ on the basis vectors by setting
\begin{equation}\label{uaction}
U_\infty x_k^\perp=h_k,\quad  U_\infty x_k=g_k\quad\text{for all $k\in \Na$},
\end{equation}
in the harmless assumption that, for each $k\in\Na$,
\begin{equation}\label{ektohnk}
 U_\infty f_k=u_k,\quad  U_\infty e_k^\perp=h_{n(k)},
\end{equation}
where $\sets{n(k)}$ is just that sequence which occurs in Condition (c).

\subsection*{Step 2.}
    The verification that $U_\infty$ in \eqref{uaction} has the desired properties
is straightforward. Fix an arbitrary $r\in\Na$ and put
$T=U_\infty S_r U_\infty ^{-1}$. Once this is done,
the index $r$ may be omitted for $S_r$, $J_r$, $Q_r$, $\varGamma_r$.

Write the Schmidt decomposition
$$
J=\sum_n s_{n}\hska{\cdot,p_{n}} q_{n},
$$
where the
$s_{n}$
are the singular values of $J$
(eigenvalues of $\left(J J^*\right)^{\frac1{2}}$),
$\sets{p_n}$, $\sets{q_n}$  are
orthonormal sets (the $p_{n}$ are eigenvectors for
$J^* J$ and the $q_{n}$ are eigenvectors for
$J J^*$).

Introduce an auxiliary operator $A$ by
\begin{equation}\label{A}
A=\sum_n s_{n}^{\frac1{4}}\hska{\cdot,p_n} q_{n},
\end{equation}
and observe that, by the Schwarz inequality,
\begin{equation}\label{schw}
\begin{gathered}
\hnorm{Af}=\hnorm{\left(J^*J\right)^{\frac18}f}\leq\hnorm{J f}^{\frac14},\\
\hnorm{A^*f}=\hnorm{\left(JJ^*\right)^{\frac18}f}\leq \hnorm{J^*f}^{\frac14}
\end{gathered}
\end{equation}
if $\norm{f}=1$.

Since $\{e_k^\perp\}_{k=1}^\infty$ is an orthonormal basis
for $(1-E)H$, \eqref{E} implies that
\begin{equation*}\label{qoper}
Q=\sum\limits_k\hska{\cdot,S^*e_k^\perp} e_k^\perp.
\end{equation*}
Whence, using  \eqref{ektohnk}, one can write
\begin{equation} \label{P}
Pf=\sum_k \ska{f,T^*h_{n(k)}} h_{n(k)}\quad (f\in L_2)
\end{equation}
where $P=U_\infty Q U_\infty ^{-1}$. By \eqref{lamgam},
\begin{equation}\label{Thk}
T^*h_{n(k)}=\sum\limits_n\hska{S^*e_k^\perp,f_n} u_n
=k\sum\limits_n\hska{e_k^\perp,\varGamma f_n} u_n \quad (k\in\Na).
\end{equation}
Prove that, for any fixed $i$, the series
\begin{equation*}\label{imtuko} 
  \sum\limits_n\hska{ e_k^\perp, \varGamma f_n} \ksk{u_n}^{(i)}(s) \quad
   (k\in\Na)
\end{equation*}
converge in the norm of $C(\R,\C)$.
Indeed, all these series are pointwise dominated on $\R$
by one series
$$
  \sum\limits_n\hnorm{\varGamma f_n}\abs{\ksk{u_n}^{(i)}(s)},
$$
which converges uniformly in $\R$ because its component subseries
\begin{gather*}
  \sum\limits_k\hnorm{\varGamma x_k}\abs{\ksk{g_k}^{(i)}(s)},\quad
  \sum\limits_k\hnorm{\varGamma x_k^\perp}\abs{\ksk{h_k}^{(i)}(s)}
\end{gather*}
are in turn dominated by the convergent series
\begin{equation*}\label{domser} 
 \sum\limits_k d(x_k)G_{k,i}, \quad
 \sum\limits_k\norm{\varGamma}H_{k,i},
\end{equation*}
respectively (see \eqref{uaction}, \eqref{defdknums}, \eqref{zndn}, \eqref{hki}).
Whence it follows via \eqref{Thk}
that, for each $k\in\Na$,
\begin{equation}\label{supest}
\norm{\ksk{T^*h_{n(k)}}^{(i)}}_{C(\R,\C)}\le C_ik,\quad
\end{equation}
with a constant $C_i$ independent of $k$.
Consider functions $\P:\R^2\to\C$, $\p:\R\to L_2$, defined, for all
$s$, $t\in\R$, by
\begin{equation} \label{Pp}
\begin{gathered}
\P(s,t)=\sum\limits_k\ksk{h_{n(k)}}(s)
         \overline{\ksk{T^*h_{n(k)}}(t)},\\
\p(s)=\overline{\P(s,\cdot)}=\sum\limits_k
\overline{\ksk{h_{n(k)}}(s)}T^*h_{n(k)}.
\end{gathered}
\end{equation}
The termwise differentiation theorem implies
that, for each $i$ and each integer $j\in[0,+\infty)$,
\begin{gather*}
\dfrac{\partial^{i+j}\P}{\partial s^i\partial t^j}(s,t)
=\sum\limits_k\ksk{h_{n(k)}}^{(i)}(s)
                             \overline{\ksk{T^*h_{n(k)}}^{(j)}(t)},\\
\dfrac{d^i\p}{ds^i}(s)=
\sum\limits_k\overline{\ksk{h_{n(k)}}^{(i)}(s)}T^*h_{n(k)},
\end{gather*}
since, by \eqref{supest} and \eqref{sumrk},
the series displayed converge (absolutely) in $C(\R^2,\C)$, $C(\R,L_2)$, respectively.
Thus, $\dfrac{\partial^{i+j}\P}{\partial s^i\partial t^j}\in C(\R^2,\C)$,
and $\dfrac{d^i\p}{ds^i}\in C(\R,L_2)$.
Observe also that, by  \eqref{sumrk}  and \eqref{Pp}, the series \eqref{P}
(viewed, of course, as one with terms belonging to $C(\R,\C)$) converges
(absolutely) in $C(\R,\C)$-norm  to the function
$$
\kskob{Pf}(s)\equiv\skoba{f,\p(s)}\equiv\int_{\R}\P(s,t)f(t)\,dt.
$$
Thus, $P$ is an integral operator, and $\P$ is its $SK^\infty$-kernel.

Since $\hnorm{S^*e_k}=\hnorm{Je_k}$ for all $k$  (see \eqref{E}),
from \eqref{null-seq} it follows
via \eqref{schw} that the operator $A$ defined in \eqref{A} is nuclear,
and hence
\begin{equation}\label{snums}
\sum_n s_{n}^{\frac1{2}}<\infty.
\end{equation}
Then, according to \eqref{A}, a kernel which induces
the nuclear operator $F=U_\infty J^* U_\infty ^{-1}$ can be represented
by the series
\begin{equation}\label{fkernels}
\sum_n s_n^{\frac12} U_\infty A^*q_n(s)
\overline{ U_\infty Ap_n(t)}
\end{equation}
convergent almost everywhere in $\R^2$. The functions used
in this bilinear expansion can be written as the series convergent in $L_2$:
$$
   U_\infty Ap_k=\sum\limits_n\hska{ p_k,A^*f_n} u_n, \quad
   U_\infty A^*q_k=\sum\limits_n\hska{q_k,Af_n} u_n\quad (k\in\Na).
$$
Show that, for any fixed $i$, the functions
$\ksk{ U_\infty Ap_k}^{(i)}$, $\ksk{ U_\infty A^*q_k}^{(i)}$ ($k\in\Na$)
make sense,
are all in $C(\R,\C)$, and their $C(\R,\C)$-norms
are bounded independent of $k$. Indeed, all the series
\begin{equation*} 
  \sum\limits_n\hska{p_k,A^*f_n}\ksk{u_n}^{(i)}(s),\quad
  \sum\limits_n\hska{q_k,Af_n}\ksk{u_n}^{(i)}(s)\quad (k\in\Na)
\end{equation*}
are dominated by one series
$$
  \sum\limits_n(\hnorm{A^*f_n}+\hnorm{Af_n})\abs{\ksk{u_n}^{(i)}(s)}.
$$
This series converges uniformly in $\R$, since it consists of two
uniformly convergent in $\R$ subseries
\begin{equation*}
\begin{gathered}
  \sum\limits_k\sk{\hnorm{A^*x_k}+\hnorm{Ax_k}}\abs{\ksk{g_k}^{(i)}(s)},\\
  \sum\limits_k\sk{\hnorm{A^*x_k^\perp}+\hnorm{Ax_k^\perp}}
  \abs{\ksk{h_k}^{(i)}(s)},
\end{gathered}
\end{equation*}
which are dominated by the following convergent series
$$
\sum\limits_k d(x_k)G_{k,i}, \quad
\sum\limits_k 2\|A\|H_{k,i},
$$
respectively (see \eqref{defdknums}, \eqref{schw}, \eqref{zndn}, \eqref{hki}).
Thus, for  functions $\F:\R^2\to \C$, $\f:\R\to L_2$, defined by
\begin{equation*}
\begin{gathered}
\F(s,t)=\sum_n s_n^{\frac12}\ksk{ U_\infty A^*q_n}(s)\overline{\ksk{ U_\infty Ap_n}(t)},\\
\f(s)=\overline{\F(s,\cdot)}=\sum_n s_n^{\frac12}\overline{\ksk{ U_\infty A^*q_n}(s)} U_\infty Ap_n, 
\end{gathered}
\end{equation*}
one can write, for all non-negative integers $i$, $j$ and all $s$, $t\in\R$,
\begin{equation*}
\begin{gathered}
\dfrac{\partial^{i+j}\F}{\partial s^i\partial t^j}(s,t)
=\sum_n s_n^{\frac12}\ksk{ U_\infty A^*q_n}^{(i)}(s)\overline{\ksk{ U_\infty Ap_n}^{(j)}(t)},\\
\dfrac{d^i\f}{ds^i}(s)
=\sum_n s_n^{\frac12}\overline{\ksk{ U_\infty A^*q_n}^{(i)}(s)} U_\infty Ap_n,
\end{gathered}
\end{equation*}
where  the series converge in $C(\R^2,\C)$, $C(\R,L_2)$, respectively,
because of \eqref{snums}. This implies that $\F$
is a $SK^\infty$-kernel of $F$.

In accordance with \eqref{splitting}, we have,
for each $f\in L_2$,
$$
\begin{gathered}
(Tf)(s)=\int_{\R} \P(s,t)f(t)\,dt+\int_{\R} \F(s,t)f(t)\,dt\\=
\int_{\R}(\P(s,t)+\F(s,t))f(t)\,dt
\end{gathered}
$$
for almost every $s$ in $\R$.
Therefore $T$ is a Carleman operator, and that kernel $\K$ of $T$, which is defined by
$\K(s,t)=\P(s,t)+\F(s,t)$ ($s$, $t\in\R$), inherits
the $SK^\infty$-kernel properties from its terms.
Consequently, $\K$ is a $SK^\infty$-kernel of $T$.

Since scalar factors do not
alter the relevant smoothness conditions, the Carleman operators
$U_\infty B_rU_\infty^{-1}=rU_\infty S_rU_\infty^{-1}$ ($r\in\Na$)
have $SK^\infty$-kernels as well.
The proof of the theorem is complete.

\section*{Acknowledgments} The author thanks the Mathematical Institute
of the  University of Jena for its hospitality, and specially
W.~Sickel and H.-J.~Schmei\ss er for useful remarks and fruitful
discussion on applying wavelets in integral representation theory.

\bibliographystyle{amsplain}

\end{document}